\documentclass[letterpaper,12pt]{article}%
\usepackage[textwidth=6.0in,textheight=9in,centering,letterpaper]{geometry}
\usepackage{amsxtra,amscd}
\usepackage{graphicx}
\usepackage[amsmath,hyperref,thmmarks]{ntheorem}
\usepackage{natbib}
\usepackage{verbatim}
\usepackage{paralist}
\usepackage{amsmath}
\usepackage{amsfonts}
\usepackage{amssymb}
\usepackage{times}
\usepackage[all,cmtip]{xy}
\usepackage{diagxy}%
\theoremnumbering{arabic}
\theoremheaderfont{\scshape}
\RequirePackage{latexsym}
\theorembodyfont{\slshape}
\theoremseparator{}
\newtheorem{X}{X}[section]

\newtheorem{corollary}[X]{Corollary}

\newtheorem{proposition}[X]{Proposition}
\newtheorem{theorem}[X]{Theorem}
\theorembodyfont{\upshape}
\newtheorem{aside}[X]{Aside}

\newtheorem{exercise}[X]{Exercise}
\newtheorem{example}[X]{Example}

\newtheorem{remark}[X]{Remark}
\newtheorem{plain}[X]{}

\theorembodyfont{\small}

\theorembodyfont{\normalsize}

\theoremstyle{nonumberplain}
\theoremheaderfont{\sc}
\theorembodyfont{\normalfont}
\theoremsymbol{\ensuremath{_\Box}}
\RequirePackage{amssymb}
\newtheorem{proof}{Proof.}
\qedsymbol{\ensuremath{_\Box}}
\theoremclass{LaTeX}
\makeindex
\newcommand{\darrow}{\!\!\downarrow\!\!}
\let\cite=\citealt

\setdefaultitem{$\diamond$\hspace{0.07in}}{}{}{}

\newcommand{\bsmall}{\begin{small}}
\newcommand{\esmall}{\end{small}}
\newcommand{\fsize}{\footnotesize}

\renewcommand{\theenumi}{{\alph{enumi}}}

\renewcommand{\theenumii}{\roman{enumii}}

\newcommand{\eb}[1]{{\itshape\bfseries#1}{\index{#1}}}
\renewcommand{\emph}{\eb}

\renewcommand{\bar}{\overline}

\renewcommand{\Gamma}{\varGamma}
\renewcommand{\Pi}{\varPi}
\renewcommand{\Sigma}{\varSigma}

\def\1{{1\mkern-7mu1}}

\DeclareMathOperator{\Aut}{Aut}
\DeclareMathOperator{\Bd}{Bd}

\DeclareMathOperator{\End}{End}

\DeclareMathOperator{\Hom}{Hom}
\DeclareMathOperator{\id}{id}

\DeclareMathOperator{\Ind}{Ind}

\DeclareMathOperator{\Ker}{Ker}

\DeclareMathOperator{\ob}{ob}

\DeclareMathOperator{\Spec}{Spec}

\newcommand{\Aff}{\mathsf{Aff}}

\newcommand{\Rep}{\mathsf{Rep}}

\newcommand{\Vc}{\mathsf{Vec}}

\newcommand{\Fib}{\textsc{Fib}}

\begin{document}

\title{Quotients of Tannakian Categories}
\author{J.S. Milne\thanks{Available at www.jmilne.org/math/.}}
\date{August 24, 2005; July 9, 2007}
\maketitle

\begin{abstract}
We classify the \textquotedblleft quotients\textquotedblright\ of a tannakian
category in which the objects of a tannakian subcategory become trivial, and
we examine the properties of such quotient categories.

\end{abstract}

\section*{Introduction}

Given a tannakian category $\mathsf{T}$ and a tannakian subcategory
$\mathsf{S}$, we ask whether there exists a quotient of $\mathsf{T}$ by
$\mathsf{S}$, by which we mean an exact tensor functor $q\colon\mathsf{T}%
\rightarrow\mathsf{Q}$ from $\mathsf{T}$ to a tannakian category $\mathsf{Q}$
such that

\begin{enumerate}
\item the objects of $\mathsf{T}$ that become trivial in $\mathsf{Q}$ (i.e.,
isomorphic to a direct sum of copies of $\1$ in $\mathsf{Q}$) are precisely
those in $\mathsf{S}$, and

\item every object of $\mathsf{Q}$ is a subquotient of an object in the image
of $q$.
\end{enumerate}

\noindent When $\mathsf{T}$ is the category $\Rep(G)$ of finite-dimensional
representations of an affine group scheme $G$ the answer is obvious: there
exists a unique normal subgroup $H$ of $G$ such that the objects of
$\mathsf{S}$ are the representations on which $H$ acts trivially, and there
exists a canonical functor $q$ satisfying (a) and (b), namely, the restriction
functor $\Rep(G)\rightarrow\Rep(H)$ corresponding to the inclusion
$H\hookrightarrow G$. By contrast, in the general case, there need not exist a
quotient, and when there does there will usually not be a canonical one. In
fact, we prove that there exists a $q$ satisfying (a) and (b) if and only if
$\mathsf{S}$ is neutral, in which case the $q$ are classified by the
$k$-valued fibre functors on $\mathsf{S}$. Here $k\overset{\text{{\tiny def}}%
}{=}\End(\1)$ is assumed to be a field.

From a slightly different perspective, one can ask the following question:
given a subgroup $H$ of the fundamental group $\pi(\mathsf{T})$ of
$\mathsf{T}$, does there exist an exact tensor functor $q\colon\mathsf{T}%
\rightarrow\mathsf{Q}$ such that the resulting homomorphism $\pi
(\mathsf{Q})\rightarrow q(\pi(\mathsf{T}))$ maps $\pi(\mathsf{Q})$
isomorphically onto $q(H)$? Again, there exists such a $q$ if and only if the
subcategory $\mathsf{T}^{H}$ of $\mathsf{T}$, whose objects are those on which
$H$ acts trivially, is neutral, in which case the functors $q$ correspond to
the $k$-valued fibre functors on $\mathsf{T}^{H}$.

The two questions are related by the \textquotedblleft tannakian
correspondence\textquotedblright\ between tannakian subcategories of
$\mathsf{T}$ and subgroups of $\pi(\mathsf{T})$ (see \ref{p8}).

In addition to proving the above results, we determine the fibre functors,
polarizations, and fundamental groups of the quotient categories $\mathsf{Q}$.

The original motivation for these investigations came from the theory of
motives (see \cite{milne2002p, milne2007rtc}).

\paragraph{Notation:}

The notation $X\approx Y$ means that $X$ and $Y$ are isomorphic, and $X\simeq
Y$ means that $X$ and $Y$ are canonically isomorphic (or that there is a given
or unique isomorphism).\qquad

\section{Preliminaries}

For tannakian categories, we use the terminology of \cite{deligneM1982}. In
particular, we write $\1$ for any identity object of a tannakian category ---
recall that it is uniquely determined up to a unique isomorphism. We fix a
field $k$ and consider only tannakian categories with $k=\End(\1)$ and only
functors of tannakian categories that are $k$-linear.

\subsection{Gerbes}

\begin{plain}
\label{p1} We refer to \cite{giraud1971}, Chapitre IV, for the theory of
gerbes. All gerbes will be for the flat (i.e., fpqc) topology on the category
$\Aff_{k}$ of affine schemes over $k$. The band (= lien) of a gerbe
$\mathcal{G}$ is denoted $\Bd(\mathcal{G})$. A commutative band can be
identified with a sheaf of groups.
\end{plain}

\begin{plain}
\label{p2} Let $\alpha\colon\mathcal{G}_{1}\rightarrow\mathcal{G}_{2}$ be a
morphism of gerbes over $\Aff_{k}$, and let $\omega_{0}$ be an object of
$\mathcal{G}_{2,k}$. Define $(\omega_{0}\darrow\mathcal{G}_{1})$ to be the
fibred category over $\Aff_{k}$ whose fibre over $S\overset{s}{\longrightarrow
}\Spec k$ has as objects the pairs $(\omega,a)$ consisting of an object
$\omega$ of $\ob(\mathcal{G}_{1,S})$ and an isomorphism $a\colon s^{\ast
}\omega_{0}\rightarrow\alpha(\omega)$ in $\mathcal{G}_{2,S}$; the morphisms
$(\omega,a)\rightarrow(\nu,b)$ are the isomorphisms $\varphi\colon
\omega\rightarrow\nu$ in $\mathcal{G}_{1,S}$ giving rise to a commutative
triangle. Thus,%
\[
\bfig\node1(0,500)[\omega]\node2(0,0)[\nu]\node3(400,250)[s^{\ast}(\omega
_{0})]\node4(1000,500)[\alpha(\omega)]\node5(1000,0)[\alpha(\nu
)]\node6(0,-200)[\mathcal{G}_{1,S}]\node7(700,-200)[\mathcal{G}_{2,S}%
]\arrow|r|[1`2;\varphi]\arrow[3`4;a]\arrow|b|[3`5;b]\arrow|r|[4`5;\alpha
(\varphi)]\arrow/{}/[5`6;{}]\arrow/{}/[6`7;{}]\efig
\]
If the map of bands defined by \noindent$\alpha$ is an epimorphism, then
$(\omega_{0}\darrow\mathcal{G}_{1})$ is a gerbe, and the sequence of bands
\begin{equation}
1\rightarrow\Bd(\omega_{0}\darrow\mathcal{G}_{1})\rightarrow\Bd(\mathcal{G}%
_{1})\rightarrow\Bd(\mathcal{G}_{2})\rightarrow1\label{e1}%
\end{equation}
\noindent is exact (\cite{giraud1971}, IV 2.5.5(i)).
\end{plain}

\begin{plain}
\label{p4} Recall (\cite{saavedra1972}, III 2.2.2) that a gerbe is said to be
tannakian if its band is locally defined by an affine group scheme. It is
clear from the exact sequence (\ref{e1}) that if $\mathcal{G}_{1}$ and
$\mathcal{G}_{2}$ are tannakian, then so also is $(\omega_{0}\darrow
\mathcal{G}_{1})$.
\end{plain}

\begin{plain}
\label{p5} The fibre functors on a tannakian category $\mathsf{T}$ form a
gerbe $\Fib(\mathsf{T})$ over $\Aff_{k}$ (\cite{deligne1990}, 1.13). Each
object $X$ of $\mathsf{T}$ defines a representation $\omega\mapsto\omega(X)$
of $\Fib(\mathsf{T})$, and in this way we get an equivalence $\mathsf{T}%
\rightarrow\Rep(\Fib(\mathsf{T}))$ of tannakian categories (\cite{deligne1989}%
, 5.11; \cite{saavedra1972}, III 3.2.3, p200). Every gerbe whose band is
tannakian arises in this way from a tannakian category (\cite[III
2.2.3]{saavedra1972}).
\end{plain}

\subsection{Fundamental groups}

\begin{plain}
\label{p6} We refer to \cite{deligne1989}, \S \S 5,6, for the theory of
algebraic geometry in a tannakian category $\mathsf{T}$ and, in particular,
for the fundamental group $\pi(\mathsf{T})$ of $\mathsf{T}$. It is the affine
group scheme\footnote{\textquotedblleft$\mathsf{T}$-sch\'{e}ma en groupes
affines\textquotedblright\ in Deligne's terminology.} in $\mathsf{T}{}$ such
that $\omega(\pi(\mathsf{T}))\simeq\underline{\Aut}^{\otimes}(\omega)$
functorially in the fibre functor $\omega$ on $\mathsf{T}$. The group
$\pi(\mathsf{T})$ acts on each object $X$ of $\mathsf{T}$, and $\omega$
transforms this action into the natural action of $\underline{\Aut}^{\otimes
}(\omega)$ on $\omega(X)$. The various realizations $\omega(\pi(\mathsf{T}))$
of $\pi(\mathsf{T})$ determine the band of $\mathsf{T}$ (i.e., the band of
$\Fib(\mathsf{T})$).
\end{plain}

\begin{plain}
\label{p7} An exact tensor functor $F\colon\mathsf{T}_{1}\rightarrow
\mathsf{T}_{2}$ of tannakian categories defines a homomorphism $\pi
(F)\colon\pi(\mathsf{T}_{2})\rightarrow F(\pi(\mathsf{T}_{1}))$
(\cite{deligne1989}, 6.4). Moreover:

\begin{enumerate}
\item $F$ induces an equivalence of $\mathsf{T}_{1}$ with a category whose
objects are the objects of $\mathsf{T}_{2}$ endowed with an action of
$F(\pi(\mathsf{T}_{1}))$ compatible with that of $\pi(\mathsf{T}_{2})$
(\cite{deligne1989}, 6.5);

\item $\pi(F)$ is flat and surjective if and only if $F$ is fully faithful and
every subobject of $F(X)$, for $X$ in $\mathsf{T}_{1}$, is isomorphic to the
image of a subobject of $X$ (cf. \cite{deligneM1982}, 2.21);

\item $\pi(F)$ is a closed immersion if and only if every object of
$\mathsf{T}_{2}$ is a subquotient of an object in the image of $q$ (ibid.).
\end{enumerate}
\end{plain}

\begin{plain}
\label{p8} For a subgroup\footnote{Note that every subgroup $H$ of
$\pi(\mathsf{T})$ is normal. For example, the fundamental group $\pi$ of the
category $\Rep(G)$ of representations of the affine group scheme $G=\Spec(A)$
is $A$ regarded as an object of $\Ind(\Rep(G))$. The action of $G$ on $A$ is
that defined by inner automorphisms. A subgroup of $\pi$ is a quotient
$A\rightarrow B$ of $A$ (as a bi-algebra) such that the action of $G$ on $A$
defines an action of $G$ on $B$. Such quotients correspond to normal subgroups
of $G$.} $H\subset\pi(\mathsf{T}){}$, we let $\mathsf{T}^{H}$ denote the full
subcategory of $\mathsf{T}$ whose objects are those on which $H$ acts
trivially. It is a tannakian subcategory of $\mathsf{T}$ (i.e., it is a
strictly full subcategory closed under the formation of subquotients, direct
sums, tensor products, and duals) and every tannakian subcategory arises in
this way from a unique subgroup of $\pi(\mathsf{T})$ (cf. \cite{bertolin2003},
1.6). The objects of $\mathsf{T}^{\pi(\mathsf{T})}$ are exactly the trivial
objects of $\mathsf{T}$, and there exists a unique (up to a unique
isomorphism) fibre functor
\[
\gamma^{\mathsf{T}}\colon\mathsf{T}^{\pi(\mathsf{T})}\rightarrow
\Vc_{k}\text{,}%
\]
namely, $\gamma^{\mathsf{T}}(X)=\Hom(\1,X).$
\end{plain}

\begin{plain}
\label{p8m}For a subgroup $H$ of $\pi(\mathsf{T})$ and an object $X$ of
$\mathsf{T}$, we let $X^{H}$ denote the largest subobject of $X$ on which the
action of $H$ is trivial. Thus $X=X^{H}$ if and only if $X$ is in
$\mathsf{T}^{H}$.
\end{plain}

\begin{plain}
\label{p9} When $H$ is contained in the centre of $\pi(\mathsf{T})$, then it
is an affine group scheme in $\mathsf{T}^{\pi(\mathsf{T})}$, and so
$\gamma^{\mathsf{T}}$ identifies it with an affine group scheme over $k$ in
the usual sense. For example, $\gamma^{\mathsf{T}}$ identifies the centre of
$\pi(\mathsf{T})$ with $\underline{\Aut}^{\otimes}(\id_{\mathsf{T}})$ (cf.
\cite{saavedra1972}, II 3.3.3.2, p150).
\end{plain}

\subsection{Morphisms of tannakian categories}

\begin{plain}
\label{p10} For a group $G$, a right $G$-object $X$, and a left $G$-object
$Y$, $X\wedge^{G}Y$ denotes the contracted product of $X$ and $Y$, i.e., the
quotient of $X\times Y$ by the diagonal action of $G$, $(x,y)g=(xg,g^{-1}y)$.
When $G\rightarrow H$ is a homomorphism of groups, $X\wedge^{G}H$ is the
$H$-object obtained from $X$ by extension of the structure group. In this last
case, if $X$ is a $G$-torsor, then $X\wedge^{G}H$ is also an $H$-torsor. See
\cite{giraud1971}, III 1.3, 1.4.
\end{plain}

\begin{plain}
\label{p11} Let $\mathsf{T}$ be a tannakian category over $k$, and assume that
the fundamental group $\pi$ of $\mathsf{T}$ is commutative. A torsor $P$ under
$\pi$ in $\mathsf{T}$ defines a tensor equivalence $\mathsf{T}{}%
\rightarrow\mathsf{T}$, $X\mapsto P\wedge^{\pi}X$, bound by the identity map
on $\Bd(\mathsf{T})$, and every such equivalence arises in this way from a
torsor under $\pi$ (cf. \cite{saavedra1972}, III 2.3). For any $k$-algebra $R$
and $R$-valued fibre functor $\omega$ on $\mathsf{T}$, $\omega(P)$ is an
$R$-torsor under $\omega(\pi)$ and $\omega(P\wedge^{\pi}X)\simeq
\omega(P)\wedge^{\omega(\pi)}\omega(X)$.
\end{plain}

\section{Quotients}

For any exact tensor functor $q\colon\mathsf{T}\rightarrow\mathsf{T}^{\prime}%
$, the full subcategory $\mathsf{T}^{q}$ of $\mathsf{T}$ whose objects become
trivial in $\mathsf{T}^{\prime}$ is a tannakian subcategory of $\mathsf{T}$ (obviously).

We say that an exact tensor functor $q\colon\mathsf{T}\rightarrow\mathsf{Q}$
of tannakian categories is a \emph{quotient functor} if every object of
$\mathsf{Q}$ is a subquotient of an object in the image of $q$; equivalently,
if the homomorphism $\pi(q)\colon\pi(\mathsf{Q})\rightarrow q(\pi\mathsf{T})$
is a closed immersion (see \ref{p7}(c)). If, in addition, the homomorphism
$\pi(q)$ is normal (i.e., its image is a normal subgroup of $q(\mathsf{T})$),
then we say that $q$ is \emph{normal}.

\begin{example}
\label{q2}Consider the exact tensor functor $\omega^{f}\colon
\Rep(G)\rightarrow\Rep(H)$ defined by a homomorphism $f\colon H\rightarrow G$
of affine group schemes. The objects of $\Rep(G)^{\omega^{f}}$ are those on
which $H$ (equivalently, the intersection of the normal subgroups of $G$
containing $f(H)$) acts trivially. The functor $\omega^{f}$ is a quotient
functor if and only if $f$ is a closed immersion, in which case it is normal
if and only if $f(H)$ is normal in $G$.
\end{example}

\begin{proposition}
\label{q1}An exact tensor functor $q\colon\mathsf{T}\rightarrow\mathsf{Q}$ of
tannakian categories is a normal quotient functor if and only if there exists
a subgroup $H$ of $\pi(\mathsf{T})$ such that $\pi(q)$ induces an isomorphism
$\pi(\mathsf{Q})\rightarrow q(H)$.
\end{proposition}

\begin{proof}
$\Longleftarrow$: Because $q$ is exact, $q(H)\rightarrow q(\pi\mathsf{T})$ is
a closed immersion. Therefore $\pi(q)$ is a closed immersion, and its image is
the normal subgroup $q(H)$ of $q(\pi\mathsf{T}).$

$\Longrightarrow$: Because $q$ is a quotient functor, $\pi(q)$ is a closed
immersion. Let $H$ be the kernel of the homomorphism $\pi(\mathsf{T}%
)\rightarrow\pi(\mathsf{T}^{q})$ defined by the inclusion $\mathsf{T}%
^{q}\hookrightarrow\mathsf{T}$. The image of $\pi(q)$ is contained in $q(H)$,
and equals it if and only if $q$ is normal. To see this, let $G=q\pi
(\mathsf{T})$, and identify $\mathsf{T}$ with the category of objects of
$\mathsf{Q}$ with an action of $G$ compatible with that of $\pi(\mathsf{Q}%
)\subset G$. Then $q$ becomes the forgetful functor, and $\mathsf{T}%
^{q}=\mathsf{T}^{\pi(\mathsf{Q})}$. Thus, $q(H)$ is the subgroup of $G$ acting
trivially on those objects on which $\pi(\mathsf{Q})$ acts trivially. It
follows that $\pi(\mathsf{Q})\subset q(H)$, with equality if and only if
$\pi(\mathsf{Q})$ is normal in $G$.
\end{proof}

In the situation of the proposition, we sometimes call $\mathsf{Q}$ \emph{a
quotient of} $\mathsf{T}$ \emph{by} $H$ (cf. \cite{milne2002p}, 1.3).

Let $q\colon\mathsf{T}\rightarrow\mathsf{Q}$ be an exact tensor functor of
tannakian categories. By definition, $q$ maps $\mathsf{T}^{q}$ into
$\mathsf{Q}^{\pi(\mathsf{Q})}$, and so we acquire a $k$-valued fibre functor
$\omega^{q}\overset{\text{{\tiny def}}}{=}\gamma^{\mathsf{Q}}\circ
(q|\mathsf{T}^{q})$ on $\mathsf{T}^{q}$:%

\[
\xymatrix{
\mathsf{T}^{q}\ar@{^{(}->}[d]\ar[r]_{q|\mathsf{T}_{q}}\ar@/^1pc/^{\omega^{q}}[rr]&
\mathsf{Q}^{\pi(\mathsf{Q})}\ar@{^{(}->}[d]\ar[r]_{\gamma^{\mathsf{Q}}}&
\Vc_{k}\\
\mathsf{T}\ar[r]^q&\mathsf{Q}\text{.}}
\]
In particular, $\mathsf{T}^{q}$ is neutral. A fibre functor $\omega$ on
$\mathsf{Q}$, defines a fibre functor $\omega\circ q$ on $\mathsf{T}$, and the
(unique) isomorphism $\gamma^{\mathsf{Q}}\rightarrow\omega|\mathsf{Q}{}%
^{\pi(\mathsf{Q})}$ defines an isomorphism $a(\omega)\colon\omega
^{q}\rightarrow(\omega\circ q)|\mathsf{T}^{q}$.

\begin{proposition}
\label{q3}Let $q\colon\mathsf{T}\rightarrow\mathsf{Q}$ be a normal quotient,
and let $H$ be the subgroup of $\pi(\mathsf{T})$ such that $\pi(\mathsf{Q}%
)\simeq q(H)$.

\begin{enumerate}
\item For $X,Y$ in $\mathsf{T}$, there is a canonical functorial isomorphism%
\[
\Hom_{\mathsf{Q}}(qX,qY)\simeq\omega^{q}(\underline{\Hom}(X,Y)^{H}).
\]

\item The map $\omega\mapsto(\omega\circ q,a(\omega))$ defines an equivalence
of gerbes
\[
r(q)\colon\Fib(\mathsf{Q})\rightarrow(\omega^{q}\darrow\Fib(\mathsf{T})).
\]

\end{enumerate}
\end{proposition}

\begin{proof}
(a) From the various definitions and \cite{deligneM1982},
\begin{align*}
\Hom_{\mathsf{Q}}(qX,qY)  &  \simeq\Hom_{\mathsf{Q}}(\1,\underline
{\Hom}(qX,qY)^{\pi(\mathsf{Q})}) &  &  \text{(ibid. 1.6.4)}\\
&  \simeq\Hom_{\mathsf{Q}}(\1,(q\underline{\Hom}(X,Y))^{q(H)}) &  &
\text{(ibid. 1.9)}\\
&  \simeq\Hom_{\mathsf{Q}}(\1,q(\underline{\Hom}(X,Y)^{H})) &  & \\
&  \simeq\omega^{q}(\underline{\Hom}(X,Y)^{H}) &  &  \text{(definition of
}\omega^{q}\text{).}%
\end{align*}

(b) The functor $\Fib(\mathsf{T})\rightarrow\Fib(\mathsf{T}^{H})$ gives rise
to an exact sequence
\[
1\rightarrow\Bd(\omega_{Q}\darrow\Fib(\mathsf{T}))\rightarrow\Bd(\mathsf{T}%
)\rightarrow\Bd(\mathsf{T}^{H})\rightarrow0
\]
(see \ref{p2}). On the other hand, we saw in the proof of (\ref{q1}) that
$H=\Ker(\pi(\mathsf{T})\rightarrow\pi(\mathsf{T}^{H}))$. On comparing these
statements, we seee that the morphism $r(q)$ of gerbes is bound by an
isomorphism of bands, which implies that it is an equivalence of gerbs (Giraud
1971, IV 2.2.6).
\end{proof}

\begin{proposition}
\label{q4}Let $(\mathsf{Q},q)$ be a normal quotient of $\mathsf{T}$. An exact
tensor functor $q^{\prime}\colon\mathsf{T}\rightarrow\mathsf{T}^{\prime}$
factors through $q$ if and only $\mathsf{T}^{q^{\prime}}\supset\mathsf{T}^{q}$
and $\omega^{q}\approx\omega^{q^{\prime}}|\mathsf{T}^{q}$.
\end{proposition}

\begin{proof}
\label{q5}The conditions are obviously necessary. For the sufficiency, choose
an isomorphism $b\colon\omega^{q}\rightarrow\omega^{q^{\prime}}|\mathsf{T}%
^{q}$. A fibre functor $\omega$ on $\mathsf{T}^{\prime}$ then defines a fibre
functor $\omega\circ q^{\prime}$ on $\mathsf{T}$ and an isomorphism
$a(\omega)|\mathsf{T}^{q}\circ b\colon\omega^{q}\rightarrow(\omega\circ
q^{\prime})|\mathsf{T}^{q}$. In this way we get a homomorphism%
\[
\Fib(\mathsf{T}^{\prime})\rightarrow(\omega^{q}\darrow\Fib(\mathsf{T}%
))\simeq\Fib(\mathsf{Q})
\]
and we can apply (\ref{p5}) to get a functor $\mathsf{Q}\rightarrow
\mathsf{T}^{\prime}$ with the correct properties.
\end{proof}

\begin{theorem}
\label{q6}Let $\mathsf{T}$ be a tannakian category over $k$, and let
$\omega_{0}$ be a $k$-valued fibre functor on $\mathsf{T}^{H}$ for some
subgroup $H\subset\pi(\mathsf{T})$. There exists a quotient $(\mathsf{Q}{},q)$
of $\mathsf{T}$ by $H$ such that $\omega^{q}\simeq\omega_{0}$.
\end{theorem}

\begin{proof}
The gerbe $(\omega_{0}\darrow\Fib(\mathsf{T}))$ is tannakian (see \ref{p4}).
From the morphism of gerbes
\[
(\omega,a)\mapsto\omega\colon(\omega_{0}\darrow\Fib(\mathsf{T}))\rightarrow
\Fib(\mathsf{T}),
\]
we obtain a morphism of tannakian categories
\[
\Rep(\Fib(\mathsf{T}))\rightarrow\Rep(\omega_{0}\darrow\Fib(\mathsf{T}))
\]
(see \ref{p5}). We define $\mathsf{Q}$ to be $\Rep(\omega_{0}%
\darrow\Fib(\mathsf{T}))$ and we define $q$ to be the composite of the above
morphism with the equivalence (see \ref{p5})
\[
\mathsf{T}\rightarrow\Rep(\Fib(\mathsf{T}))\text{.}%
\]
Since a gerbe and its tannakian category of representations have the same
band, an argument as in the proof of Proposition \ref{q3} shows that $\pi(q)$
maps $\pi(\mathsf{Q{}})\mathsf{\ }$isomorphically onto $q(H)$. A direct
calculation shows that $\omega^{q}$ is canonically isomorphic to $\omega_{0}$.
\end{proof}

We sometimes write $\mathsf{T}/\omega$ for the quotient of $\mathsf{T}$
defined by a $k$-valued fibre functor $\omega$ on a subcategory of
$\mathsf{T}$.

\begin{example}
\label{q7}Let $(\mathsf{T},w,\mathbb{T}{})$ be a Tate triple, and let
$\mathsf{S}$ be the full subcategory of $\mathsf{T\ }$of objects isomorphic to
a direct sum of integer tensor powers of the Tate object $\mathbb{T}$. Define
$\omega_{0}$ to be the fibre functor on $\mathsf{S}$,
\[
X\mapsto\varinjlim_{n}\Hom(\bigoplus_{-n\leq r\leq n}\1(r),X).
\]
Then the quotient tannakian category $\mathsf{T}/\omega_{0}$ is that defined
in \cite{deligneM1982}, 5.8.
\end{example}

\begin{remark}
\label{q7m}Let $q\colon\mathsf{T}\rightarrow\mathsf{Q}$ be a normal quotient
functor. Then $\mathsf{T}$ can be recovered from $\mathsf{Q}$, the
homomorphism $\pi(\mathsf{Q})\rightarrow q(\pi(\mathsf{T}))$, and the actions
of $q(\pi(\mathsf{T}))$ on the objects of $\mathsf{Q}$ (apply \ref{p7}(a)).
\end{remark}

\begin{remark}
\label{q7n}A fixed $k$-valued fibre functor on a tannakian category
$\mathsf{T}$ determines a Galois correspondence between the subsets of
$\ob(\mathsf{T})$ and the equivalence classes of quotient functors
$\mathsf{T}\rightarrow\mathsf{Q}$.
\end{remark}

\begin{exercise}
\label{q7p}Use (\ref{p10}, \ref{p11}) to express the correspondence between
fibre functors on tannakian subcategories of $\mathsf{T}$ and normal quotients
of $\mathsf{T}$ in the language of $2$-categories.
\end{exercise}

\begin{aside}
\label{q8}Let $G$ be the fundamental group $\pi(\mathsf{T})$ of a tannakian
category $\mathsf{T}$, and let $H$ be a subgroup of $G$. We use the same
letter to denote an affine group scheme in $\mathsf{T}$ and the band it
defines. Then, under certain hypotheses, for example, if all the groups are
commutative, there will be an exact sequence%
\[
\cdots\rightarrow H^{1}(k,G)\rightarrow H^{1}(k,G/H)\rightarrow H^{2}%
(k,H)\rightarrow H^{2}(k,G)\rightarrow H^{2}(k,G/H).
\]
The category $\mathsf{T}$ defines a class $c(\mathsf{T})$ in $H^{2}(k,G)$,
namely, the $G$-equivalence class of the gerbe of fibre functors on
$\mathsf{T}$, and the image of $c(\mathsf{T})$ in $H^{2}(k,G/H)$ is the class
of $\mathsf{T}^{H}$. Any quotient of $\mathsf{T}$ by $H$ defines a class in
$H^{2}(k,H)$ mapping to $c(\mathsf{T})$ in $H^{2}(k,G)$. Thus, the exact
sequence suggests that a quotient of $\mathsf{T}$ by $H$ will exist if and
only if the cohomology class of $\mathsf{T}^{H}$ is neutral, i.e., if and only
if $\mathsf{T}^{H}$ is neutral as a tannakian category, in which case the
quotients are classified by the elements of $H^{1}(k,G/H)$ (modulo
$H^{1}(k,G)$). When $\mathsf{T}$ is neutral and we fix a $k$-valued fibre
functor on it, then the elements of $H^{1}(k,G/H)$ classify the $k$-valued
fibre functors on $\mathsf{T}^{H}$. Thus, the cohomology theory suggests the
above results, and in the next subsection we prove that a little more of this
heuristic picture is correct.
\end{aside}

\subsection{The cohomology class of the quotient}

For an affine group scheme $G$ over a field $k$, $H^{r}(k,G)$ denotes the
cohomology group computed with respect to the flat topology. When $G$ is not
commutative, this is defined only for $r=0,1,2$ (\cite{giraud1971}).

\begin{proposition}
\label{q9}Let $(\mathsf{Q},q)$ be a quotient of $\mathsf{T}$ by a subgroup $H$
of the centre of $\pi(\mathsf{T})$. Suppose that $\mathsf{T}$ is neutral, with
$k$-valued fibre functor $\omega$. Let $G=\underline{\Aut}^{\otimes}(\omega)$,
and let $\wp(\omega^{q})$ be the $G/\omega(H)$-torsor $\underline{\Hom}%
(\omega|\mathsf{T}^{H},\omega^{q})$. Under the connecting homomorphism
\[
H^{1}(k,G/H)\rightarrow H^{2}(k,H)
\]
the class of $\wp(\omega^{q})$ in $H^{1}(k,G/H)$ maps to the class of
$\mathsf{Q}$ in $H^{2}(k,H)$.
\end{proposition}

\begin{proof}
Note that $H=\Bd(\mathsf{Q})$, and so the statement makes sense. According to
\cite{giraud1971}, IV 4.2.2, the connecting homomorphism sends the class of
$\wp(\omega^{q})$ to the class of the gerbe of liftings of $\wp(\omega^{q})$,
which can be identified with $(\omega^{q}\darrow\Fib(\mathsf{T}))$. Now
Proposition \ref{q3} shows that the $H$-equivalence class of $(\omega
^{q}\darrow\Fib(\mathsf{T}))$ equals that of $\Fib(\mathsf{Q})$ which (by
definition) is the cohomology class of $\mathsf{Q}$.
\end{proof}

\subsection{Semisimple normal quotients}

Everything can be made more explicit when the categories are semisimple.
Throughout this subsection, $k$ has characteristic zero.

\begin{proposition}
\label{q10}Every normal quotient of a semisimple tannakian category is semisimple.
\end{proposition}

\begin{proof}
A tannakian category is semisimple if and only if the identity component of
its fundamental group is pro-reductive (cf. \cite{deligneM1982}, 2.28), and a
normal subgroup of a reductive group is reductive (because its unipotent
radical is a characteristic subgroup).
\end{proof}

Let $\mathsf{T}$ be a semisimple tannakian category over $k$, and let
$\omega_{0}$ be a $k$-valued fibre functor on a tannakian subcategory
$\mathsf{S}$ of $\mathsf{T}.$ We can construct an explicit quotient
$\mathsf{T}/\omega_{0}$ as follows. First, let $(\mathsf{T}/\omega
_{0})^{\prime}$ be the category with one object $\bar{X}$ for each object $X$
of $\mathsf{T}$, and with
\[
\Hom_{(\mathsf{T}/\omega_{0})^{\prime}}(\bar{X},\bar{Y})=\omega_{0}%
(\underline{\Hom}(\bar{X},\bar{Y})^{H})
\]
where $H$ is the subgroup of $\pi(\mathsf{T})$ defining $\mathsf{S}$. There is
a unique structure of a $k$-linear tensor category on $(\mathsf{T}/\omega
_{0})^{\prime}$ for which $q\colon\mathsf{T}\rightarrow(\mathsf{T}{}%
/\omega_{0})^{\prime}$ is a tensor functor. With this structure,
$(\mathsf{T}/\omega_{0})^{\prime}$ is rigid, and we define $\mathsf{T}%
{}/\omega_{0}$ to be its pseudo-abelian hull. Thus, $\mathsf{T}/\omega_{0}$
has%
\[%
\begin{array}
[c]{rl}%
\text{objects:} & \text{pairs }(\bar{X},e)\text{ with }X\in\ob(\mathsf{T}%
{})\text{ and }e\text{ an idempotent in }\End(\bar{X})\text{,}\\
\text{morphisms:} & \Hom_{\mathsf{T}/\omega_{0}}((\bar{X},e),(\bar
{Y},f))=f\circ\Hom_{(\mathsf{T}/\omega_{0})^{\prime}}(\bar{X},\bar{Y})\circ e.
\end{array}
\]
Then $(\mathsf{T}/\omega_{0},q)$ is a quotient of $\mathsf{T}$ by $H$, and
$\omega^{q}\simeq\omega_{0}$.

Let $\omega$ be a fibre functor on $\mathsf{T}$, and let $a$ be an isomorphism
$\omega_{0}\rightarrow\omega|\mathsf{T}^{H}$. The pair $(\omega,a)$ defines a
fibre functor $\omega_{a}$ on $\mathsf{T}/\omega_{0}$ whose action on objects
is determined by%
\[
\omega_{a}(\bar{X})=\omega(X)
\]
and whose action on morphisms is determined by
\[
\bfig
\node a(0,500)[\Hom(\bar{X},\bar{Y})]
\node b(2000,500)[\Hom(\omega_{a}(\bar{X}),\omega_{a}(\bar{Y}))]
\node c(0,0)[\omega_{0}(\underline{\Hom}(X,Y)^{H})]
\node d(1000,0)[\omega(\underline{\Hom}(X,Y)^{H})]
\node e(2000,0)[\underline{\Hom}(\omega(X),\omega(Y))^{\omega(H)}]
\arrow/-->/[a`b;\omega_{a}]
\arrow|r|/=/[a`c;\text{def}]
\arrow[c`d;a]
\arrow[d`e;\simeq]
\arrow/^{(}->/[e`b;{}]
\efig
\]

\noindent The map $(\omega,a)\mapsto\omega_{a}$ defines an equivalence
$(\omega_{0}\darrow\Fib(\mathsf{T}))\rightarrow\Fib(\mathsf{T}/\omega_{0})$.

Let $H_{1}\subset H_{0}\subset\pi(\mathsf{T})$, and let $\omega_{0}$ and
$\omega_{1}$ be $k$-valued fibre functors on $\mathsf{T}^{H_{0}}$ and
$\mathsf{T}{}^{H_{1}}$ respectively. A morphism $\alpha\colon\omega
_{0}\rightarrow\omega_{1}|\mathsf{T}^{H_{0}}$ defines an exact tensor functor
$\mathsf{T}{}/\omega_{0}\rightarrow\mathsf{T}/\omega_{1}$ whose action on
objects is determined by
\[
\bar{X}\text{ (in }\mathsf{T}^{H_{0}}\text{) }\mapsto\bar{X}\text{ (in
}\mathsf{T}^{H_{1}}\text{),}%
\]
and whose action on morphisms is determined by
\[
\bfig
\node a(0,500)[\Hom_{\mathsf{T}/\omega_{0}}(\bar{X},\bar{Y})]
\node b(2000,500)[\Hom_{\mathsf{T}/\omega_{1}}(\bar{X},\bar{Y})]
\node c(0,0)[\omega_{0}(\underline{\Hom}_{\mathsf{T}}(X,Y)^{H_{0}})]
\node d(1000,0)[\omega_{1}(\underline{\Hom}_{\mathsf{T}}(X,Y)^{H_{0}})]
\node e(2000,0)[\omega_{1}(\underline{\Hom}_{\mathsf{T}}(X,Y)^{H_{1}}))]
\arrow/-->/[a`b;{}]
\arrow|r|/=/[a`c;\text{def}]
\arrow[c`d;\alpha]
\arrow/^{(}->/[d`e;{}]
\arrow/=/[e`b;\text{def}]
\efig
\]
When $H_{1}=H_{0}$, this is an isomorphism (!) of tensor categories
$\mathsf{T}/\omega_{0}\rightarrow\mathsf{T}/\omega_{1}$.

Let $(\mathsf{Q}_{1},q_{1})$ and $(\mathsf{Q}_{2},q_{2})$ be quotients of
$\mathsf{T}$ by $H$. For simplicity, assume that $\pi\overset
{\text{{\tiny def}}}{=}\pi(\mathsf{T})$ is commutative. Then $\underline
{\Hom}(\omega^{q_{1}},\omega^{q_{2}})$ is $\pi/H$-torsor, and we assume that
it lifts to a $\pi$-torsor $P$ in $\mathsf{T}$, so $P\wedge^{\pi}%
(\pi/H)=\underline{\Hom}(\omega^{q_{1}},\omega^{q_{2}})$. Then%

\[
\mathsf{T}\xrightarrow{X\mapsto P\wedge^{\pi}X}\mathsf{T}%
\xrightarrow{q_{2}}\mathsf{Q}_{2}%
\]
realizes $\mathsf{Q}_{2}$ as a quotient of $\mathsf{T}$ by $H$, and the
corresponding fibre functor on $\mathsf{T}^{H}$ is $P\wedge^{\pi}\omega
^{q_{2}}\simeq\omega^{q_{1}}$. Therefore, there exists a commutative diagram
of exact tensor functors%
\[
\begin{CD} \mathsf{T}@>{X\mapsto P\wedge ^{\pi}X}>>\mathsf{T} \\
@VV{q_{1}}V @VV{q_{2}}V \\
\mathsf{Q}_{1} @>\phantom{X\mapsto P\wedge ^{\pi}X}>>\mathsf{Q}_{2},\end{CD}
\]
which depends on the choice of $P$ lifting $\underline{\Hom}(\omega^{q_{1}%
},\omega^{q_{2}})$ in an obvious way.

\section{Polarizations}

We refer to \cite{deligneM1982}, 5.12, for the notion of a (graded)
polarization on a Tate triple over $\mathbb{R}{}$. We write $\mathsf{V}$ for
the category of $\mathbb{Z}{}$-graded complex vector spaces endowed with a
semilinear automorphism $a$ such that $a^{2}v=(-1)^{n}v$ if $v\in V^{n}$. It
has a natural structure of a Tate triple (ibid. 5.3). The canonical
polarization on $\mathsf{V}$ is denoted $\Pi^{\mathsf{V}}$.

A morphism $F\colon(\mathsf{T}_{1},w_{1},\mathbb{T}_{1})\rightarrow
(\mathsf{T}_{2},w_{2},\mathbb{T}_{2})$ of Tate triples is an exact tensor
functor $F\colon\mathsf{T}_{1}\rightarrow\mathsf{T}_{2}$ preserving the
gradations together with an isomorphism $F(\mathbb{T}_{1})\simeq
\mathbb{\mathbb{T}}_{2}$. We say that such a morphism is \emph{compatible
}with graded polarizations $\Pi_{1}$ and $\Pi_{2}$ on $\mathsf{T}_{1}$ and
$\mathsf{T}_{2}$ (denoted $F\colon\Pi_{1}\mapsto\Pi_{2}$) if
\[
\psi\in\Pi_{1}(X)\Rightarrow F\psi\in\Pi_{2}(FX)\text{,}%
\]
in which case, for any $X$ homogeneous of weight $n$, $\Pi_{1}(X)$ consists of
the sesquilinear forms $\psi\colon X\otimes\bar{X}\rightarrow\1(-n)$ such that
$F\psi\in\Pi_{2}(FX)$. In particular, given $F$ and $\Pi_{2}$, there exists at
most one graded polarization $\Pi_{1}$ on $\mathsf{T}_{1}$ such that
$F\colon\Pi_{1}\mapsto$ $\Pi_{2}$.

Let $\mathsf{T}=(\mathsf{T},w,\mathbb{T)}$ be an algebraic Tate triple over
$\mathbb{R}{}$ such that $w(-1)\neq1$. Given a graded polarization $\Pi$ on
$\mathsf{T}$, there exists a morphism of Tate triples $\xi_{\Pi}%
\colon\mathsf{T}\rightarrow\mathsf{V}$ (well defined up to isomorphism) such
that $\xi_{\Pi}\colon\Pi\mapsto\Pi^{\mathsf{V}}$ (\cite{deligneM1982}, 5.20).
Let $\omega_{\Pi}$ be the composite%
\[
\mathsf{T}^{w(\mathbb{G}_{m})}\overset{\xi_{\Pi}}{\rightarrow}\mathsf{V}%
^{w(\mathbb{G}_{m})}\overset{\gamma^{\mathsf{V}}}{\rightarrow}\Vc_{\mathbb{R}%
{}};
\]
it is a fibre functor on $\mathsf{T}^{w(\mathbb{G}_{m})}$.

\subsubsection{A criterion for the existence of a polarization}

\begin{proposition}
\label{z1}Let $\mathsf{T=(T}{},w,\mathbb{T)}$ be an algebraic Tate triple over
$\mathbb{R}{}$ such that $w(-1)\neq1$, and let $\xi\colon\mathsf{T}%
\rightarrow\mathsf{V}$ be a morphism of Tate triples. There exists a graded
polarization $\Pi$ on $\mathsf{T}$ (necessarily unique) such that $\xi
\colon\Pi\mapsto\Pi^{\mathsf{V}}$ if and only if the real algebraic group
$\underline{\Aut}^{\otimes}(\gamma^{\mathsf{V}}\circ\xi|\mathsf{T}%
{}^{w(\mathbb{G}_{m})})$ is anisotropic.
\end{proposition}

\begin{proof}
Let $G=\underline{\Aut}^{\otimes}(\gamma^{\mathsf{V}}\circ\xi|\mathsf{T}%
{}^{w(\mathbb{G}_{m})})$. Assume $\Pi$ exists. The restriction of $\Pi$ to
$\mathsf{T}^{w(\mathbb{G}_{m})}$ is a symmetric polarization, which the fibre
functor $\gamma^{\mathsf{V}}\circ\xi$ maps to the canonical polarization on
$\Vc_{\mathbb{R}{}}$. This implies that $G$ is anisotropic (\cite{deligne1972}%
, 2.6).

For the converse, let $X$ be an object of weight $n$ in $\mathsf{T}%
{}_{(\mathbb{C}{})}$. A sesquilinear form $\psi\colon\xi(X)\otimes
\overline{\xi(X)}\rightarrow\1(-n)$ arises from a sesquilinear form on $X$ if
and only if it is fixed by $G$. Because $G$ is anisotropic, there exists a
$\psi\in\Pi^{\mathsf{V}}(\xi(X))$ fixed by $G$ (ibid., 2.6), and we define
$\Pi(X)$ to consist of all sesquilinear forms $\phi$ on $X$ such that
$\xi(\phi)\in\Pi^{\mathsf{V}}(\xi(X))$. It is now straightforward to check
that $X\mapsto\Pi(X)$ is a polarization on $\mathsf{T}$.
\end{proof}

\begin{corollary}
\label{z2}Let $F\colon(\mathsf{T}_{1},w_{1},\mathbb{T}_{1})\rightarrow
(\mathsf{T}_{2},w_{2},\mathbb{T}_{2})$ be a morphism of Tate triples, and let
$\Pi_{2}$ be a graded polarization on $\mathsf{T}_{2}$. There exists a graded
polarization $\Pi_{1}$ on $\mathsf{T}_{1}$ such that $F\colon\Pi_{1}\mapsto
\Pi_{2}$ if and only if the real algebraic group $\underline{\Aut}^{\otimes
}(\gamma^{\mathsf{V}}\circ\xi_{\Pi_{2}}\circ F|\mathsf{T}_{1}^{w(\mathbb{G}%
_{m})})$ is anisotropic.
\end{corollary}

\subsubsection{Polarizations on quotients}

The next proposition gives a criterion for a polarization on a Tate triple to
pass to a quotient Tate triple.

\begin{proposition}
\label{z3}Let $\mathsf{T}$ $=(\mathsf{T},w,\mathbb{T)}$ be an algebraic Tate
triple over $\mathbb{R}{}$ such that $w(-1)\neq1$. Let $(\mathsf{Q,}q)$ be a
quotient of $\mathsf{T}{}$ by $H\subset\pi(\mathsf{T})$, and let $\omega^{q}$
be the corresponding fibre functor on $\mathsf{T}^{H}$. Assume $H\supset
w(\mathbb{G}_{m})$, so that $\mathsf{Q}$ inherits a Tate triple structure from
that on $\mathsf{T}$, and that $\mathsf{Q}$ is semisimple. Given a graded
polarization $\Pi$ on $\mathsf{T}$, there exists a graded polarization
$\Pi^{\prime}$ on $\mathsf{Q}$ such that $q\colon\Pi\mapsto\Pi^{\prime}$ if
and only if $\omega^{q}\approx\omega_{\Pi}|\mathsf{T}^{H}$.
\end{proposition}

\begin{proof}
$\Rightarrow$: Let $\Pi^{\prime}$ be such a polarization on $\mathsf{Q}$, and
consider the functors
\[
\mathsf{T}\overset{q}{\rightarrow}\mathsf{Q}\overset{\xi_{\Pi^{\prime}}%
}{\rightarrow}\mathsf{V},\quad\xi_{\Pi^{\prime}}\colon\Pi^{\prime}\mapsto
\Pi^{\mathsf{V}}.
\]
\noindent Both $\xi_{\Pi^{\prime}}\circ q$ and $\xi_{\Pi}$ are compatible with
$\Pi$ and $\Pi^{\mathsf{V}}$ and with the Tate triple structures on
$\mathsf{T}$ and $\mathsf{V}$, and so $\xi_{\Pi^{\prime}}\circ q\approx
\xi_{\Pi}$ (\cite{deligneM1982}, 5.20). On restricting everything to
$\mathsf{T}^{w(\mathbb{G}_{m})}$ and composing with $\gamma^{V},$ we get an
isomorphism $\omega_{\Pi^{\prime}}\circ(q|\mathsf{T}^{w(\mathbb{G}_{m}%
)})\approx\omega_{\Pi}$. Now restrict this to $\mathsf{T}^{H}$, and note that
\[
\left(  \omega_{\Pi^{\prime}}\circ(q|\mathsf{T}^{w(\mathbb{G}_{m})})\right)
|\mathsf{T}{}^{H}=(\omega_{\Pi^{\prime}}|\mathsf{Q}^{\pi(\mathsf{Q})}%
)\circ(q|\mathsf{T}^{H})\simeq\omega^{q}%
\]
because $\omega_{\Pi^{\prime}}|\mathsf{Q}^{\pi(\mathsf{Q})}\simeq\gamma^{Q}$.

$\Leftarrow$: The choice of an isomorphism $\omega^{q}\rightarrow\omega_{\Pi
}|\mathsf{T}^{H}$ determines an exact tensor functor
\[
\mathsf{T}/\omega^{q}\rightarrow\mathsf{T}/\omega_{\Pi}\text{.}%
\]
\noindent As the quotients $\mathsf{T}/\omega^{q}$ and $\mathsf{T}/\omega
_{\Pi}$ are tensor equivalent respectively to $\mathsf{Q}$ and $\mathsf{V}$,
this shows that there is an exact tensor functor $\xi\colon\mathsf{Q}%
{}\rightarrow\mathsf{V}$ such that $\xi\circ q\approx\xi_{\Pi}$. Evidently
$\underline{\Aut}^{\otimes}(\gamma^{\mathsf{V}}\circ\xi|\mathsf{Q}%
^{w(\mathbb{G}_{m})})$ is isomorphic to a subgroup of $\underline
{\Aut}^{\otimes}(\gamma^{\mathsf{V}}\circ\xi_{\Pi}|\mathsf{T}^{w(\mathbb{G}%
_{m})})$. Since the latter is anisotropic, so also is the former
(\cite{deligne1972}, 2.5). Hence $\xi$ defines a graded polarization
$\Pi^{\prime}$ on $\mathsf{Q}$ (Proposition \ref{z1}), and clearly $q\colon
\Pi\mapsto\Pi^{\prime}$.
\end{proof}

\bsmall {
\bibliographystyle{cbe}
\bibliography{T:/MData/bib/refs}

\begin{thebibliography}{}

\bibitem[\protect\astroncite{Bertolin}{2003}]{bertolin2003}
{\sc Bertolin, C.} 2003.
\newblock Motivic {G}alois theory for motives of niveau $\leq$ 1.
\newblock Preprint arXiv:math.NT/0309379.

\bibitem[\protect\astroncite{Deligne}{1972}]{deligne1972}
{\sc Deligne, P.} 1972.
\newblock La conjecture de {W}eil pour les surfaces {$K3$}.
\newblock {\em Invent. Math.} 15:206--226.

\bibitem[\protect\astroncite{Deligne}{1989}]{deligne1989}
{\sc Deligne, P.} 1989.
\newblock Le groupe fondamental de la droite projective moins trois points, pp.
  79--297.
\newblock {\em In} Galois groups over ${\bf Q}$ (Berkeley, CA, 1987), volume~16
  of {\em Math. Sci. Res. Inst. Publ.} Springer, New York.

\bibitem[\protect\astroncite{Deligne}{1990}]{deligne1990}
{\sc Deligne, P.} 1990.
\newblock Cat{\'e}gories tannakiennes, pp. 111--195.
\newblock {\em In} The Grothendieck Festschrift, Vol.\ {II}, Progr. Math.
  Birkh{\"a}user Boston, Boston, MA.

\bibitem[\protect\astroncite{Deligne and Milne}{1982}]{deligneM1982}
{\sc Deligne, P. and Milne, J.~S.} 1982.
\newblock Tannakian categories, pp. 101--228.
\newblock {\em In} Hodge cycles, motives, and {S}himura varieties, Lecture
  Notes in Mathematics. Springer-Verlag, Berlin.

\bibitem[\protect\astroncite{Giraud}{1971}]{giraud1971}
{\sc Giraud, J.} 1971.
\newblock Cohomologie non ab{\'e}lienne.
\newblock Springer-Verlag, Berlin.

\bibitem[\protect\astroncite{Milne}{2002}]{milne2002p}
{\sc Milne, J.~S.} 2002.
\newblock Polarizations and {G}rothendieck's standard conjectures.
\newblock {\em Ann. of Math. (2)} 155:599--610.

\bibitem[\protect\astroncite{Milne}{2007}]{milne2007rtc}
{\sc Milne, J.~S.} 2007.
\newblock Rational {T}ate classes on abelian varieties.
\newblock In preparation.

\bibitem[\protect\astroncite{Saavedra~Rivano}{1972}]{saavedra1972}
{\sc Saavedra~Rivano, N.} 1972.
\newblock Cat{\'e}gories {T}annakiennes.
\newblock Springer-Verlag, Berlin.

\end{thebibliography}
} \esmall

\fsize

\bigskip\noindent J.S. Milne, Ann Arbor, MI, USA

\end{document}